\documentclass[12pt]{article}

\usepackage{amssymb,amsmath,exscale}

\usepackage{graphicx}

\usepackage[applemac]{inputenc}
\usepackage[T1]{fontenc}

\usepackage{color}
\definecolor{marin}{rgb}   {0.,   0.3,   0.7} 
\definecolor{rouge}{rgb}   {0.8,   0.,   0.} 
\definecolor{sepia}{rgb}   {0.8,   0.5,   0.} 
\usepackage[colorlinks,citecolor=marin,linkcolor=rouge,
            bookmarksopen,
            bookmarksnumbered
           ]{hyperref}

\newtheorem{lemma}{Lemma}[section]

\newtheorem{theorem}[lemma]{Theorem}
\newtheorem{proposition}[lemma]{Proposition}
\newtheorem{corollary}[lemma]{Corollary}
\newtheorem{remark}[lemma]{Remark}
\newtheorem{example}[lemma]{Example}

\newtheorem{notation}[lemma]{Notation}
\newtheorem{definition}[lemma]{Definition}
\newtheorem{conclusion}[lemma]{Conclusion}
\numberwithin{equation}{section}
\newcommand{\QED}{\mbox{}\hfill \raisebox{-0.2pt}{\rule{5.6pt}{6pt}\rule{0pt}{0pt}} 
          \medskip\par}             
\newenvironment{Proof}{\noindent
    \parindent=0pt\abovedisplayskip = 0.5\abovedisplayskip
    \belowdisplayskip=\abovedisplayskip{\bfseries Proof. }}{\QED}

\newcommand{\ad}{\mathrm{ad}}

\newcommand{\dd}{\mathrm{d}}

\newcommand{\N}{\mathbb{N}}

\newcommand{\R}{\mathbb{R}}

\newcommand{\C}{\mathbb{C}}

\newcommand{\T}{\mathbb{T}}
\newcommand{\Z}{\mathbb{Z}}

\newcommand{\Norm}[2]{\|#1\|\left.\vphantom{T_{j_0}^0}\!\!\right._{#2}}         
             
\title{Modified energy for split-step methods applied to the linear Schr\"odinger equation }        
\author{Arnaud Debussche and Erwan Faou\\[4ex]
\small INRIA \& ENS Cachan Bretagne,  Avenue Robert Schumann
F-35170 Bruz\\
\small\it email: \tt Arnaud.Debussche@bretagne.ens-cachan.fr, Erwan.Faou@inria.fr\\[4ex]
}      
\begin{document}
\maketitle
\abstract{
We consider the linear Schr\"odinger equation and its discretization by split-step methods where the part corresponding to the Laplace operator is approximated by the midpoint rule. We show that the numerical solution coincides with the exact solution of a modified partial differential equation at each time step. This shows the existence of a  modified energy preserved by the numerical scheme. This energy is close to the exact energy if the numerical solution is smooth.  As a consequence, we give uniform regularity estimates for the numerical solution over arbitrary long time. \\[2ex]
{\bf MSC numbers}: 65P10, 37M15}\\[2ex]
{\bf Keywords}: Schr\"odinger equation, Splitting integrators, Long-time behavior, Backward error analysis.

%newpage

%\tableofcontents

%%%%%%%%%%%%%%%%%%%%%%%%%%%%%%
\section{Introduction}
%%%%%%%%%%%%%%%%%%%%%%%%%%%%%%

We consider the linear Schr\"odinger equation
\begin{equation}
\label{E1}
\partial_t u(t,x) = -i \Delta u (t,x) + iV(x) u(t,x),\quad u(0,x) = u^0(x),
\end{equation}
with initial condition $u^0$, and potential function $V(x) \in \R$.
The wave function $u(x,t)$ depends on  $x \in \T^d$ or $\R^d$ and the time $t > 0$. The operator $\Delta$ is the $d$-dimensional Laplace operator. In the following, we consider mainly the case where $x \in \T^d$. The case of the whole space is totally similar. 
The equation \eqref{E1} is symplectic and its solution preserves the $L^2$ norm and the energy
\begin{equation}
\label{Eorig}
u \mapsto \int_{\T^d} |\nabla u|^2 + V |u|^2 \dd x = \langle u|-\Delta+V|u\rangle.  
\end{equation}
The solution of \eqref{E1} is given by 
$$
u(t,x) = \exp(i t (-\Delta + V)) u^0(x),
$$
and a standard method to simulate this solution is to consider the approximation
\begin{equation}
\label{EsplitD}
\exp(i h (-\Delta + V)) \simeq \exp(-ih \Delta) \exp(ih V)
\end{equation}
for a small stepsize $h > 0$. 
The solution at a given time $t = nh$ is then approximated by 
\begin{equation}
\label{Eapprox}
 \exp(i t (-\Delta + V)) u^0 \simeq \Big(\exp(-ih \Delta) \exp(ih V)\Big)^n u^0. 
\end{equation}
The advantage of this method is that it yields a symplectic scheme preserving the $L^2$ norm. Moreover, it is very easy to implement by using the fast Fourier transform: while the operator $\Delta$ is diagonal in the Fourier space, the operator $V$ acts as a multiplication operator in the phase space.  
For finite time, this splitting scheme yields a consistent numerical scheme: as $h \to 0$ and if the numerical solution is smooth, it can be shown that \eqref{Eapprox} yields a convergent approximation of order $1$ in $h$, see \cite{Jahn00}. Considering higher order approximation such as the symmetric Strang splitting or higher order splitting methods allows to obtain higher order approximation scheme under the assumption that the numerical solution is smooth enough, see \cite{Jahn00,Hans08}. 

Concerning the long-time behaviour of such methods, very few results exist. In \cite{DF07}, {\sc Dujardin \& Faou} showed the conservation of the regularity of the numerical solution \eqref{Eapprox} in $\T^1$ over very long time, provided the potential function is small and smooth. Moreover, even in this situation, resonances effects appear for some values of $h$: typically when $\exp(-ih \Delta)$ posseses eigenvalues close to $1$. 

In the finite dimensional case, the long time behaviour of splitting method can be understood upon using the Baker-Campbell-Hausdorff formula (see for instance \cite{HLW}). Roughly speaking, this result states that for two matrices $A$ and $B$, we can write
$$
\exp(t A ) \exp(tB) = \exp( t Z(t))
$$
where $Z(t) = A + B + t [A,B] + t^2 \cdots$, with $[A,B] = AB - BA$ the matrix commutator. Hence the long time behaviour of the numerical solution corresponding to \eqref{Eapprox} can be analyzed by considering the properties of the matrix $Z(t)$ which is a small perturbation of the original operator $A +B$ for small time $t$. However, to be valid, the BCH formula requires $h$ to be small enough with respect to the inverse of the norms of $A$ and $B$. This makes this strategy impossible to apply directly for unbounded operators, unless a drastic CFL like condition is used for the full discretization of \eqref{E1}. 

In this paper, we consider the time discretization 
\begin{equation}
\label{Esplit}
\exp(ih (-\Delta + V)) \simeq  \exp(ihV) R(-ih\Delta) 
\end{equation}
where 
$$
R(z) = \frac{1+z/2}{1-z/2}
$$
is the stability function of the midpoint rule. Such an approximation is clearly consistent with \eqref{E1} if the solution is smooth enough. Moreover, it defines a symplectic numerical scheme preserving the $L^2$ norm, and easily implemented by using the fast Fourier transform. Similar schemes have been considered in \cite{Ascher, Stern, Zhang}. 

Recall that for all $x\in \R$ we have 
$$
\frac{1 + ix }{1 - ix} = \exp(2i\arctan(x)). 
$$
and hence we can write
$$
R(-ih\Delta) = \frac{1 -ih\Delta/2}{1 + ih\Delta/2} = \exp(2 i \arctan\big(-\frac{h\Delta}{2}\big)),  
$$
where now $ 2 \arctan\big(-\frac{h\Delta}{2}\big)$ is a bounded operator from $L^2$ to itself. Using this representation, we show in this work that there exists a symmetric operator $S(h): L^2 \to L^2$ such that 
$$
 \exp(ihV) R(-ih\Delta) = \exp(ih S(h)),
$$
with
$$
S(h) = -\frac{2}{h}\arctan\big(\frac{h\Delta}{2}\big) + \tilde{V}(h)
$$
where $\tilde{V}(h): L^2 \to L^2$ is a modified potential. 

Hence, for all $n$ and all initial value $ u^0$, we have 
$$
 u^n = \big(\exp(ihV) R(-ih\Delta)\big)^n u^0 = \exp(inh S(h)) u^0
$$
and hence the numerical solution $u^n$ coincides with the exact solution of the {\em modified equation}
$$
\partial_t u = S(h) u
$$
at each time step $t_n = nh$. This implies that the associated energy
 $$
\langle u \, | \, S(h) \, | \, u\rangle
$$
is preserved along the numerical solution associated with the split-step scheme \eqref{Esplit}. Moreover this energy is close to the original energy \eqref{Eorig} if $u$ is smooth. Using these properties, we give  regularity bounds for the numerical solution over arbitrary long time. 

Such a result is to our knowledge the first extension in an infinite dimensional setting of the classical backward error analysis for Hamiltonian ordinary differential equation (see \cite{HLW,Reic04}). Note in particular that as in the case of {\em linear} ordinary differential equation, this result is valid for arbitrary long time, while such results classically hold for exponentially long time with respect to the step size for nonlinear ordinary differential equations. 

It is worth noticing that such result does not hold hold for the splitting scheme \eqref{EsplitD} for which it is known that resonance effects occur, see \cite{DF07}.  The main difference between \eqref{Esplit} and \eqref{EsplitD} lies in the high frequencies regularization effect of the midpoint rule: by essence, the logarithm of the operator $R(-ih\Delta)$ is bounded while the logarithm of $\exp(-ih \Delta)$ is not well defined when $h \Delta$ possesses eigenvalues close to multiples of $2\pi$. Note that this does not affect the approximation property of the scheme for finite time and smooth numerical solution. 

Similarly this result does not automatically extend to situations where the propagator $R(-ih\Delta)$ is replaced by a higher order approximation of $\exp(-ih\Delta)$, or for higher order splitting schemes (see \cite[Chap III]{HLW}). We discuss this point in the last section of this work, and show by numerical experiments that in general resonance effects appear. 

Let us mention that in the nonlinear situation, results exist concerning the long-time behaviour of splitting scheme applied to the nonlinear Schr\"odinger equation: see the recent works of \textsc{Faou, Gr\'ebert \& Paturel} \cite{FGP1,FGP2} and \textsc{Gauckler \& Lubich} \cite{GL08a,GL08b} for the long time behaviour of splitting schemes applied to NLS when the initial solution is small. However,  to our knowledge no existence results for a global modified energy have been proved. Note that in this direction, concerning the numerical approximation of solitary wave,  \textsc{Duran \& Sanz-Serna} \cite{Duran00} have proved the existence of a modified solitary wave over finite time for the numerical solution associated with the midpoint rule. 

%%%%%%%%%%%%%%%%%%%%%%%%%%%%%%%%%%
\section{Statement of the results}
%%%%%%%%%%%%%%%%%%%%%%%%%%%%%%%%%%

We represent a function $u \in L^2(\T^d)$ by its Fourier coefficients $u = (u_k)_{k \in \Z^d}$ defined as 
$$
u_k = \frac{1}{(2\pi)^d}\int_{\T^d} u(x) e^{i k \cdot x} \dd x
$$
where for $k = (k_1,\ldots,k_d) \in \Z^d$ and $x = (x_1,\cdots,x_d) \in \T^d$ we set $k \cdot x = k_1 x_1 + \cdots k_d x_d$. 
We define
$$
\Norm{u}{}^2 = \sum_{k \in \Z^d} |u_k|^2,\quad \mbox{and}\quad \Norm{u}{H^s}^2 = \sum_{k \in \Z^d} (1 + |k|^2)^s |u_k|^2
$$
the $L^2$ and the $H^s$ Sobolev norms on $\T^d$, where for $k = (k_1,\ldots,k_d) \in \Z^d$, we set 
$$
|k|^2 = k_1^2 + \cdots k_d^2. 
$$
For an operator $A = (A_{k\ell})_{k,\ell \in \Z^d}$ acting in the Fourier space $\C^{\Z^d}$ and for $\alpha > 1$ we set 
$$
\Norm{A}{\alpha} = \sup_{k,\ell} |A_{k\ell}| \big(1+  |k - \ell|^\alpha\big). 
$$
We denote by 
$$
\mathcal{L}_{\alpha} = \{ A = (A_{k\ell})_{k,\ell \in \Z^d}\, | \, \Norm{A}{\alpha} < \infty\, \}. 
$$
If $A \in \mathcal{L}_{\alpha}$ with $\alpha > d$, we can easily show that $A \in \mathcal{L}(L^2)$: see Lemma \ref{ELB} below. 

We say that $A$ is symmetric if for all $k,\ell \in \Z^d$, we have $A_{k\ell} = \overline{A}_{\ell k}$, or equivalently $A^* = A$. In this situation, for $u \in L^2$, we set 
$$
\langle u | \, A \, | u \rangle = \sum_{k,\ell \in \Z^d} \bar{u}_k A_{k\ell} u_\ell = (u,Au) \in \R
$$
where $(\,\cdot \, , \, \cdot\, )$ is the $L^2$ product in $\T^d$. 
For two operators $A$ and $B$, we set 
$$
\ad_A(B) = AB - BA. 
$$
Finally, with a real function $W(x)$ we associate the operator $W = (W_{k\ell})_{k,\ell \in \Z^d}$ with components $W_{k\ell}= W_{k- \ell}$ where $W_n$ denote the Fourier coefficient of $W$ associated with  $n \in \Z^d$. Thus the operator $(W_{k\ell})_{k,\ell \in \Z^d}$ acting in the Fourier space corresponds to the multiplication by $W$. Note moreover that with this identification, $\Norm{W}{\alpha} < \infty$ with $\alpha > d$ implies that $\Norm{W}{L^{\infty}} < \infty$. 

The goal of this paper is to prove the following results: 

\begin{theorem}
\label{T1}
Let $ \alpha > d$, and 
assume that $\Norm{V}{\alpha} < \infty$. 
There exist $h_0> 0$ and a constant $C$ such that for all $h \in (0, h_0)$, 
there exists a symmetric operator $S(h)$ such that 
$$
\exp(ih V) R(-ih\Delta) = \exp (ih S(h)), 
$$
satisfying for all $h$, 
$$
S(h) = -\frac{2}{h} \arctan\big(\frac{h \Delta}{2}\big) + V(h) + h W(h)
$$
where $V(h)$ and $W(h)$ satisfy
\begin{equation}
\label{truc1}
\Norm{V(h)}{\alpha} + \Norm{W(h)}{\alpha} \leq C \Norm{V}{\alpha},
\end{equation}
 and 
where moreover $V(h)$ is given by the convergent series in $\mathcal{L}_{\alpha}$
\begin{equation}
\label{truc2}
V(h) = \big(\dd \exp_{Z_0(h)}\big)^{-1} (V) = V + \sum_{k \geq 1} \frac{B_k}{k!} i^k \ad_{Z_0(h)}^k (V)
\end{equation}
with $Z_0(h) = -2  \arctan\displaystyle\big(\frac{h \Delta}{2}\big) $, and where the $B_k$ are the Bernouilli numbers. 
\end{theorem}

\begin{remark}
The size of $h_0$ is only proportional to the inverse of $\Norm{V}{\alpha}$, and hence is a reasonably small parameter. In particular it does not depend on a possible space discretization of the problem through a CFL condition. 
\end{remark}

The following result shows that $S(h)$ defines a ``modified'' energy when applied to smooth functions: 
\begin{proposition}
\label{P1}
Let $\beta \in [0,1]$. 
Assume that $u \in H^{1 +\beta}(\T^d)$, then we have for $h \in (0, h_0)$, 
\begin{equation}
\label{Ecomp}
\big|\langle u | S(h) | u \rangle -  \langle u | -\Delta + V | u \rangle \big|\leq C h^\beta \Norm{u}{H^{1 + \beta}}^2. 
\end{equation}
where $C$ depends on $\beta$ and $V$. 
\end{proposition}

The next results shows the conservation the modified energy $S(h)$ along the numerical solution associated with the split-step propagator. As a consequence, we give a regularity bound for the numerical solution over arbitrary long time.  
\begin{corollary}
\label{C1}
Assume that $u^0 \in L^2(\T^d)$ and $h \in (0,h_0)$. For all $n \geq 1$, we define
$$
u^n = \big(\exp(ih V) R(-ih\Delta)\big)^n u^0. 
$$
Then for all $n$ we have
\begin{equation}
\label{Eenerg}
\langle u^n | S(h) | u^n \rangle = \langle u^0 | S(h) | u^0 \rangle. 
\end{equation}
If moreover $u^0 \in H^1$, 
then there exists a constant $C_0$ depending on  $V$ and $\alpha$ such that for all $n \in \N$,  
\begin{equation}
\label{Ebreg}
\sum_{ |k| \leq 1/\sqrt{h}} |k|^2|u^n_k|^2 + \frac{1}{h}  \sum_{ |k| > 1/\sqrt{h}} |u^n_k|^2 \leq C_0 \Norm{u^0}{H^1}^2. 
\end{equation}
\end{corollary}

This last result shows that $H^1$ estimate are preserved over arbitrary long time only for ``low'' modes $|k| < 1/\sqrt{h}$ whereas the remaining high frequencies part is small in $L^2$.

\begin{remark}
The results above obviously remain valid when considering the full discretization of \eqref{E1} by collocation methods (see for instance \cite{L08}), with estimates independent of the spectral discretization parameter. 
\end{remark}

\begin{remark}
\label{NONO}
The previous results easily extend to the splitting scheme
$$
R(-ih \Delta) \exp(ihV)
$$
and to the Strang splitting
\begin{equation}
\label{Estrang}
\exp(ihV/2) R(-ih \Delta) \exp(ihV/2). 
\end{equation}
Note that in this last situation, the fact that the method is of order $2$ allows to take $\beta \in [0,2]$ in \eqref{Ecomp}. See Section \ref{Section7} for further details on other possible extensions. 
\end{remark}

%%%%%%%%%%%%%%%%%%%%%
\section{Formal series}
%%%%%%%%%%%%%%%%%%%%%

We now start the proof of Theorem \ref{T1}. 

In the following, we set 
$$
Z_0 := -2 \arctan\big(\frac{h\Delta}{2}
\big)
$$
the diagonal operator with coefficients
$$
\lambda_k = (Z_0)_{kk} = 2 \arctan\big(\frac{h|k|^2}{2}\big), \quad k \in \Z^d. 
$$
We look for a function $ t \to Z(t)$ taking value into the set of operator acting on $\C^{\Z^d}$ such that $Z(0) = Z_0$ and 
$$
\forall\, t \in [0,h],\quad  e^{itV} e^{iZ_0} = e^{iZ(t)}. 
$$ 
Derivating the equation in $t$, this yields (see \cite{HLW})
$$
i V e^{it V }e^{iZ_0}= i \big(\dd \exp_{iZ(t)} Z'(t) \big) e^{iZ(t)}. 
$$
Hence $Z(t)$ has to satisfy the equation (see \cite[Chap. III.4]{HLW})
\begin{equation}
\label{EZt}
Z'(t) =  (\dd \exp_{iZ(t)})^{-1} V  =  i \sum_{k\geq 0} \frac{B_k}{k!} \mathrm{ad}_{iZ(t)}^k(V).
\end{equation}
and $Z(0) =  Z_0$. Here, the $B_k$ are the Bernouilli numbers. Recall that for $z \in \C$, $|z| < 2\pi$, the expression 
$$
\sum_{k \geq 0} \frac{B_k}{k!} z^k = \frac{z}{e^{z} - 1}
$$
defines a power series of radius $2\pi$. 

We define the formal series 
$$
Z (t) = \sum_{\ell \geq 0 } t^\ell Z_\ell
$$
where $Z_\ell$, $\ell \geq 1$, are unknown operators. 

Plugging this expression into \eqref{EZt} we find
$$
\begin{array}{rcl}
\displaystyle\sum_{\ell \geq 1} \ell t^{\ell - 1} Z_{\ell} 
&=&  \displaystyle \sum_{k \geq 0} \frac{B_k}{k!} \Big( i \sum_{\ell \geq 0} t^\ell \ad_{Z_\ell} \Big)^k(V)\\[3ex]
&=&  \displaystyle \sum_{\ell \geq 0} t^\ell \sum_{k \geq 0} \frac{B_k}{k!} i^k  \sum_{\ell_1 + \cdots + \ell_k = \ell } \ad_{Z_{\ell_1}} \cdots \ad_{Z_{\ell_k}}(V).
\end{array}
$$
Identifying the coefficients in the formal series, we find the induction formula: 
\begin{equation}
\label{Erec}
\forall\, \ell \geq 1,\quad
(\ell+1) Z_{\ell+1} =  \sum_{k \geq 0} \frac{B_k}{k!} i^k  \sum_{\ell_1 + \cdots + \ell_k = \ell } \ad_{Z_{\ell_1}} \cdots \ad_{Z_{\ell_k}}(V). 
\end{equation}
Note that we easily show by induction that for all $\ell$,  $Z_\ell$ is symmetric. For $\ell = 1$, this equation yields
\begin{equation}
\label{EZ1}
Z_1 =  \sum_{k \geq 0} \frac{B_k}{k!} i^k \mathrm{ad}_{Z_0}^k(V).
\end{equation}
Note that the main difference with the finite dimensional situation is that the ``first'' term in the expansion is given by an infinite series and that it depends on the small parameter $h$ through the operator $Z_0$. The key to control this term is to estimate the norm of the operator $\ad_{Z_0}$. 

%%%%%%%%%%%%%%%%%%%%%%%%%%%%%
\section{Proof of Theorem \ref{T1}}
%%%%%%%%%%%%%%%%%%%%%%%%%%%%%

\begin{lemma}
\label{Lprod}
Assume that $\alpha > d$. 
There exist a constant $C_\alpha$ such that  for all operator $A$ and $B$, $$
\Norm{AB}{\alpha} \leq C_\alpha\Norm{A}{\alpha}\Norm{B}{\alpha}. 
$$
\end{lemma}
\begin{Proof}
We have for $k,\ell \in \Z^d$, 
$$
\begin{array}{rcl}
|(AB)_{k\ell} |( 1+| k - \ell|^\alpha) &\leq& (1+ | k - \ell|^\alpha) \displaystyle\sum_{p\in \Z^d} |A_{kp}| |B_{kp}|\\[2ex]
&\leq &\displaystyle \Norm{A}{\alpha}\Norm{B}{\alpha}  \sum_{p\in \Z^d} \frac{1+ | k - \ell|^\alpha}{(1+ |k-p|^\alpha)(1+ |p - \ell|^\alpha)}
\end{array}
$$
But as the function $x \to x^{\alpha}$ is convex for $x > 0$, we have 
$$
1 + |k-p|^\alpha \leq 1 + \big(|k-\ell| + |\ell - p|\big)^\alpha \leq 2^{\alpha - 1} \big(1 + |k-\ell|^\alpha + 1+ |\ell - p|^\alpha\big). 
$$
Hence we have 
$$
|(AB)_{k\ell} |(1+ | k - \ell|^\alpha) \leq 2^{\alpha-1} \Norm{A}{\alpha}\Norm{B}{\alpha}  \sum_{p\in \Z^d}\Big( \frac{1}{1+ |k-p|^\alpha} + \frac{1}{ 1+ |p - \ell|^\alpha}\Big)
$$
and this shows the result, as $\alpha > d$. 
\end{Proof}

\begin{lemma}
\label{ELB}
Let $\alpha> d$. There exist a constant $M_\alpha$ such that for all symmetric operator $B$ and for all $u \in L^2$, we have 
$$
|\langle u  |  B | u \rangle | \leq M_\alpha \Norm{B}{\alpha} \Norm{u}{}^2. 
$$
\end{lemma}
\begin{Proof}
We have 
$$
\begin{array}{rcl}
|\langle u  |  B |\, u \rangle | &\leq& \displaystyle\sum_{k,\ell} |B_{k\ell}| |u_k||u_\ell|\\[2ex]
&\leq & \Norm{B}{\alpha}\displaystyle \sum_{k,\ell} \frac{1}{1+ |k - \ell|^\alpha} |u_k||u_\ell|\\[2ex]
&\leq&  \Norm{B}{\alpha} \displaystyle \sum_{k,\ell} \frac{1}{1+ |k - \ell|^\alpha} |u_k|^2 
\end{array}
$$
using the formula $|u_k||u_\ell| \leq \frac{1}{2}(|u_k|^2 + |u_\ell|^2) $. 
This yields the result. 
\end{Proof}

\begin{lemma}
\label{L1}
Recall that $Z_0 = \displaystyle 2 \arctan\big(\frac{h\Delta}{2}\big)$, and 
let $W = (W_{k\ell})_{k,\ell \in \Z^d}$ be an operator. We have for all $\alpha > 1$
\begin{equation}
\label{Enono}
\Norm{\ad_{Z_0} W}{\alpha} \leq \pi \Norm{W}{\alpha}. 
\end{equation}
\end{lemma}
\begin{Proof}
For $k,\ell \in \Z^d$ we have as $Z_0$ is diagonal
$$
\begin{array}{rcl}
\big(\mathrm{ad}_{Z_0}W \big)_{k\ell} &=&  (\lambda_k - \lambda_\ell)W_{k\ell}, \\[2ex]
&=&  \big(2\arctan(h |k|^2/2) - 2\arctan(h|\ell|^2/2)\big)W_{k\ell}. 
\end{array}
$$
Hence we have for all $k,\ell \in \Z^d$, 
$$
\left|\big(\mathrm{ad}_{Z_0}W \big)_{k\ell} \right| \leq \pi |W_{k\ell}|
$$
and this shows the result. 
\end{Proof}

Using this Lemma, we see using \eqref{EZ1} that 
\begin{equation}
\label{truc4}
\Norm{Z_1}{\alpha} \leq \Norm{V}{\alpha} \sum_{k \geq 0} \frac{|B_k|}{k!} \pi^k \leq C \Norm{V}{\alpha}
\end{equation}
is bounded. 
In components, we calculate using the expression of $\ad_{Z_0}$ that
\begin{equation}
\label{EZ1nono}
(Z_1)_{k\ell} =  V_{k\ell} \frac{i(\lambda_k - \lambda_\ell)}{\exp(i(\lambda_k - \lambda_\ell)) - 1}
\end{equation}

Note that for any bounded operator $A$ and $B$, we always have 
$$
\Norm{\ad_A(B)}{\alpha}  \leq 2C_\alpha \Norm{A}{\alpha}\Norm{B}{\alpha} 
$$
where $C_\alpha$ is given by Lemma \ref{Lprod}
We define now the following numbers: 
$$
\zeta_0 = \pi\quad \mbox{and}\quad \zeta_\ell = 2 C_\alpha \Norm{Z_\ell}{\alpha},\quad \mbox{for}\quad \ell \geq 1. 
$$
Using \eqref{Erec} and Lemma \ref{L1}, we easily see that we have the estimates
$$
\forall\, \ell\geq 1,\quad \frac{1}{2C_\alpha}(\ell+1) \zeta_{\ell+1} \leq  \Norm{V}{\alpha}\sum_{k \geq 0} \frac{|B_k|}{k!}  \sum_{\ell_1 + \cdots + \ell_k = \ell } \zeta_{\ell_1}\cdots \zeta_{\ell_k} . 
$$
Now for any $\rho$ such that $\pi < \rho < 2\pi$, there exist a constant $M$ such that for all $k$, $|B_k| \leq k! M \rho^{-k}$.  Hence we can write
$$
\forall\, \ell\geq 1,\quad \frac{1}{2C_\alpha}(\ell+1) \zeta_{\ell+1} \leq  M \Norm{V}{\alpha}\sum_{k \geq 0}  \rho^{-k }\sum_{\ell_1 + \cdots + \ell_k = \ell } \zeta_{\ell_1}\cdots \zeta_{\ell_k} . 
$$
Let $\zeta(t)$ be the formal series $\zeta(t) = \sum_{\ell \geq 0} t^\ell \zeta_\ell$. 
Multiplying the previous equation by $t^\ell$ and summing over $\ell \geq 0$, we find 
$$
\frac{1}{2C_\alpha}\zeta'(t) \leq M \Norm{V}{\alpha} \sum_{k \geq 0}  \rho^{-k}\zeta(t)^k = M \Norm{V}{\alpha} \frac{1}{1 - \zeta(t)/\rho}. 
$$
Let $\eta(t)$ be the solution of the differential equation: 
$$
\eta'(t) = 2 M C_\alpha \Norm{V}{\alpha} \frac{1}{1 - \eta(t)/\rho}, \quad \eta(0) = \pi. 
$$
Taking $\rho = 3\pi /2$, we easily see that for $t \leq \frac{\pi}{32 M C_\alpha\Norm{V}{\alpha}}$, the solution can be written 
$$
\eta(t) = \frac{3\pi}{2} \left( 1 - \sqrt{\frac19 - \frac{16}{3} M C_\alpha\Norm{V}{\alpha} t}\right),
$$
and defines an analytic function of $t$. Expanding $\eta(t) = \sum_{\ell \geq 0} t^\ell\eta_\ell$, we see that the coefficients satisfy the relations $\eta_0 = \pi$ and 
$$
\forall\, \ell\geq 1,\quad \frac{1}{2C_\alpha}(\ell+1) \eta_{\ell+1} =  M \Norm{V}{}\sum_{k \geq 0}  \rho^{-k }\sum_{\ell_1 + \cdots + \ell_k = \ell } \eta_{\ell_1}\cdots \eta_{\ell_k} 
$$
with $\rho = \frac{3\pi}{2}$. 
By induction, this shows that $\zeta_\ell \leq \eta_\ell$. Moreover, for all $z \in \C$ with $|z| \leq \frac{\pi}{32 M C_\alpha \Norm{V}{\alpha}}$, we have as the coefficients $\zeta_\ell$ are positive, 
$$
|\zeta(z)| = \left|\sum_{\ell = 0}^\infty \zeta_\ell z^\ell \right| \leq \sum_{\ell = 0}^\infty \zeta_\ell |z|^\ell = \zeta(|z|) \leq \eta(|z|) \leq \frac{3\pi}{2}. 
$$
Using Cauchy estimates, we see that 
$$
\forall\, \ell \geq 1, \quad \Norm{Z_\ell}{} = \frac{1}{2C_\alpha} \zeta_\ell = \frac1{2C_\alpha} \frac{\zeta^{(\ell)}(0)}{\ell! } \leq \frac{3\pi}{4C_\alpha}  \Big(\frac{32 M C_\alpha\Norm{V}{\alpha}}{\pi} \Big)^\ell. 
$$

The theorem is now proved by setting 
$$
V(h) = Z_1,\quad \mbox{and}\quad W(h) = \sum_{\ell \geq 2} h^{\ell - 2} Z_\ell
$$
which defines a convergent power series for $|h| < h_0 = \frac{\pi}{32 M C_\alpha \Norm{V}{\alpha}}$. The estimate \eqref{truc1} on $V(h)$ is then an easy consequence of \eqref{truc4}. The estimate \eqref{truc1} on $W(h)$ is easily proved. 

%%%%%%%%%%%%%%%%%%%%%%%%
\section{Modified energy}
%%%%%%%%%%%%%%%%%%%%%%%%%

We give now the proof of Proposition \ref{P1}. 

For all $x \in \R$, we have
$$
 \arctan(x) -  x = -\int_{0}^x \frac{y^2}{1 + y^2} \dd y. 
$$
For $k \in \Z^d$, this yields
$$
\frac{2}{h}\arctan\big(\frac{h|k|^2}{2}\big) - |k|^2 = - \frac{2}{h} \int_{0}^{h|k|^2/2} \frac{y^2}{1 + y^2} \dd y. 
$$
Let $\gamma \in [0,2]$, it is clear that for all $y \in \R$, 
$$
 \frac{y^2}{1 + y^2} \leq y^{\gamma}. 
$$
Hence we have for all $k \in \Z^d$, 
$$
\Big|\frac{2}{h}\arctan\big(\frac{h|k|^2}{2}\big) - |k|^2 \Big|\leq  \frac{2}{h} \int_{0}^{h|k|^2/2} y^{\gamma}\dd y \leq C h^{\gamma} |k|^{2\gamma + 2}. 
$$
This shows that for all $v$, 
\begin{equation}
\label{Ebdelta}
\Big| \langle v |  -\frac{2}{h} \arctan\big(\frac{h\Delta}{2}\big) |  v \rangle - \langle v  |  -\Delta  |  v \rangle\Big|
\leq C h^{\gamma} \Norm{v}{H^{1 + \gamma}}^2. 
\end{equation}
Now we have 
$$
 \langle v \, | \, V(h) \, | \, v \rangle - \langle v \, | \, V \, | \, v \rangle =  \displaystyle\sum_{k \geq 1 } \frac{B_k}{k!}
 \langle v \, | \, i^k \ad_{Z_0(h)}^k(V)  \, | \, v \rangle
$$
Recall that $Z_0(h) = - 2 \arctan\big(\frac{h\Delta}{2}\big)$ is a positive operator. The operator $Z_0(h)^{1/2}$ is hence well defined, and for an operator $W$ we have in components
$$
(Z_0(h)^{1/2}W)_{k\ell} = \Big(2 \arctan\big(\frac{h|k|^2}{2}\big)\Big)^{1/2} W_{k\ell}. 
$$
Hence we have for all $\alpha > 1$, 
$$
\Norm{Z_0(h)^{1/2}W}{\alpha} \leq \sqrt{\pi} \Norm{W}{\alpha}
\quad \mbox{and}\quad\Norm{W Z_0(h)^{1/2}}{\alpha} \leq \sqrt{\pi} \Norm{W}{\alpha}. 
$$
Now using Lemma \ref{ELB} and the fact that $Z_0(h)$ is symmetric,  we have for all $v$ and all operator $W$
\begin{multline*}
| \langle v \, | \,\ad_{Z_0(h)}(W) \, | v \rangle | \leq (\Norm{Z_0(h)^{1/2}W}{\alpha} +\Norm{W Z_0(h)^{1/2}}{\alpha}) \Norm{Z_0(h)^{1/2} v}{} \Norm{v}{}\\[2ex]
\leq  2 \sqrt{\pi}\Norm{W}{\alpha} \Norm{Z_0(h)^{1/2} v}{} \Norm{v}{}. 
\end{multline*}
Hence we have 
\begin{multline*}
\big| \langle v \, | \, V(h) \, | \, v \rangle - \langle v \, | \, V \, | \, v \rangle\big| \leq 
\displaystyle 2 \sum_{k \geq 1 } \frac{|B_k|}{k!} \pi^{k-1/2} \Norm{V}{\alpha}
\Norm{Z_0(h)^{1/2} v}{} \Norm{v}{}\\[2ex]
\leq C \Norm{V}{\alpha}\Norm{Z_0(h)^{1/2} v}{} \Norm{v}{}
\end{multline*} 
Using \eqref{Ebdelta} with $\gamma = 0$, this shows that 
$$
\big| \langle v \, | \, V(h) \, | \, v \rangle - \langle v \, | \, V \, | \, v \rangle\big|  \leq C \Norm{V}{\alpha} h \Norm{u}{H^1}\Norm{u}{}. 
$$
Finally, we easily have using \eqref{truc1} that
$$
\big| \langle v \, | \, W(h) \, | \, v \rangle \big| \leq C\Norm{V}{\alpha} h \Norm{u}{}^2. 
$$
Summing the previous inequalities with $\gamma = \beta$ in \eqref{Ebdelta}  we have that 
$$
\langle\, u | S(h) | u \rangle -  \langle\, u | \Delta + V | u \rangle \leq C h^\beta \Norm{u}{H^{1 + \beta}}^2 + C \Norm{V}{\alpha} h \Norm{u}{H^1}\Norm{u}{}
$$
and this yields the result.

%%%%%%%%%%%%%%%%%%%%%%%%%%%%%%%%%%
\section{Bounds for the numerical solution}
%%%%%%%%%%%%%%%%%%%%%%%%%%%%%%%%%%

We prove now Corollary \ref{C1}. Note that Eqn. \eqref{Eenerg} is classic. 

Using the fact that $V$ is symmetric, we have for all $n$, $\Norm{u^n}{} = \Norm{u^0}{}$ where $\Norm{\cdot}{}$ denotes the $L^2$ norm. 

Using Lemma \ref{ELB}, we can write for all $v \in L^2$, 
$$
 \langle v  |  S(h)  |  v \rangle  =   \frac{1}{h} \langle v  |  -2 \arctan\big(\frac{h\Delta}{2}\big) |  v \rangle +  \langle v  |  V(h) + h W(h) \, | \, v \rangle   
$$
whence using \eqref{truc1}, Lemma \ref{ELB} and the fact that $Z_0$ is a positive operator
$$
 |\langle v  | \, S(h)  |  v \rangle|  \geq \frac{1}{h} \langle v  |  -2 \arctan\big(\frac{h\Delta}{2}\big) |  v \rangle - C \Norm{V}{\alpha} \Norm{v}{}^2. 
$$
Hence using \eqref{Eenerg} we have that for all $n$, 
$$
\begin{array}{rcl}
\displaystyle\frac{1}{h} \langle u^n | - 2 \arctan\big(\frac{h\Delta}{2}\big)|  u^n \rangle  &\leq&  \langle u^n  |  S(h)  |  u^n\rangle + C \Norm{V}{\alpha} \Norm{u^n}{}^2\\[2ex]
&\leq &  \langle u^0  |  S(h) |  u^0\rangle + C \Norm{V}{\alpha} \Norm{u^0}{}^2. 
\end{array}
$$
Using \eqref{Ecomp} with $\beta = 0$, we find that there exists a constant such that for all $n$, 
\begin{equation}
\label{Eatan}
\frac{1}{h} \langle u^n  | - 2 \arctan\big(\frac{h\Delta}{2}\big) |  u^n \rangle \leq C_0\Norm{u^0}{H^1}^2. 
\end{equation}
Now we have for all $x> 0$
\begin{equation}
\label{Ebornarc}
x > \frac12 \Longrightarrow\arctan x > \arctan\big(\frac12\big)\quad \mbox{and}\quad x \leq \frac12 \Longrightarrow \arctan x > \frac{2x}{3}. 
\end{equation}
Applying this inequality to \eqref{Eatan} by considering the set of frequencies $h|k|^2 \leq 1$ and $h|k|^2 > 1$ immediately yields the result. 

\section{Higher order approximations}

\label{Section7}

In this section we further investigate the long time behaviour by numerical simulations and consider higher-order numerical schemes. 

We perform the simulations with $d = 1$, $u^0 = 2/(2 - \cos(x))$ and $V(x) = \cos(x) + \sin(6x)$. In the next figures, we show the maximal size of the oscillations of the truncated $H^1$ norm 
\begin{equation}
\label{EH120}
\Big(\sum_{k = -20}^{20} (1 + |k|^2) |u_k^n|^2 \Big)^{1/2}
\end{equation}
along the numerical solution $u^n$ from $t = 0$ to $t = 50$, and for stepsize ranging from $h = 0.01$ to $h = 0.1$. 

As expected, we see that this quantity is uniformly bounded for the splitting scheme \eqref{Esplit} (Figure 1). 
\begin{figure}[ht]
\begin{center}
\rotatebox{0}{\resizebox{0.8\linewidth}{0.2\linewidth}{%
   \includegraphics{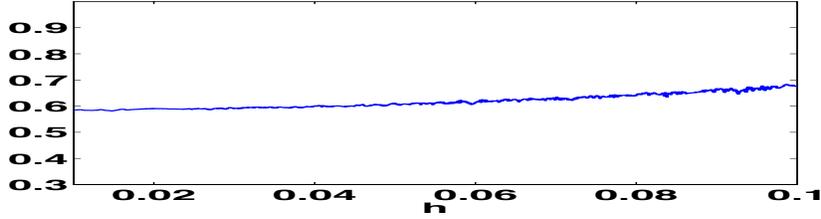}}}
   %\rotatebox{0}{\resizebox{!}{0.3\linewidth}{%
   %\includegraphics{fig8.eps}}}
   \end{center}
   \caption{Midpoint approximation of the exponential.}
\end{figure}

As explained in Remark \ref{NONO}, our methods easily extends  to the Strang splitting scheme \eqref{Estrang}. 
Considering the alternative Strang splitting
$$
R(-ih\Delta/2) \exp(-ihV) R(-ih\Delta/2),
$$
the same argument does not apply straightforwardly. The obstruction occurs in Lemma \ref{L1} where $R(-ih\Delta)$ is replaced by $R(-ih\Delta/2)^2$ in the definition of the operator $Z_0$,  transforming $\pi$ by $2\pi$ in inequality \eqref{Enono}. 

Nevertheless, as shown in Figure 2, the same uniform conservation phenomenon can be observed. This might be justified using the fact that the operator $Z_1$ defined in \eqref{EZ1nono} still makes sense in this situation.

Next we consider schemes of the form 
\begin{equation}
\label{EsplitN}
\exp(ihV) \prod_{j = 1}^s R(-\gamma_j h \Delta)
\end{equation}
where $\gamma_j \in \R$, $j = 1,\ldots,s$ are coefficients satisfying $\gamma_1 + \ldots + \gamma_s = 1$. Such an approximation will be a higher order approximation of the splitting scheme \eqref{EsplitD} for suitable $\gamma_j$ satisfying given algebraic conditions (see for instance \cite[Chap III]{HLW}). Of course, all these schemes remain symplectic and preserve the $L^2$ norm.

\begin{figure}[ht]
\begin{center}
\rotatebox{0}{\resizebox{0.8\linewidth}{0.2\linewidth}{%
   \includegraphics{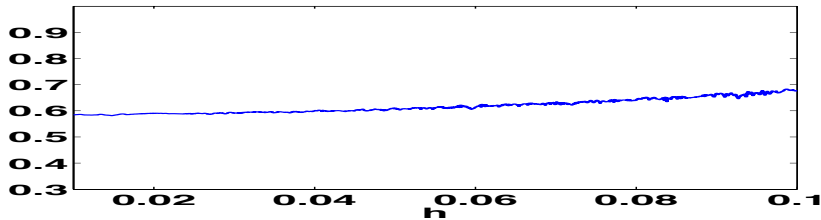}}}
   %\rotatebox{0}{\resizebox{!}{0.3\linewidth}{%
   %\includegraphics{fig8.eps}}}
   \end{center}
   \caption{Strang splitting $R(-ih\Delta/2) \exp(-ihV) R(-ih\Delta/2)$.}
\end{figure}

In Figures 3, 4 and 5, we consider successively classical symmetric composition methods of order $4$, $6$ and $8$ (see \cite[Chap V]{HLW} and the references therein). The method of order $4$ is the triple jump method for which $s= 3$, 
\begin{equation}
\label{Egammas}
\gamma_1 = \gamma_3 = \frac{1}{2 - 2^{1/3}}, \quad \mbox{and}\quad \gamma_2 = -\frac{2^{1/3}}{2 - 2^{1/3}}. 
\end{equation}
The methods of order $6$ corresponds to the methods given by {\sc Yoshida} (see \cite{Yoshida} and \cite[Section V.3.2]{HLW}) and requires $s = 7$, while the method of order $8$ is the methods given by {\sc Suzuki \& Umeno}, see \cite{Suzuki}, and requires $s = 15$. 

What we observe is that for the method of order $4$, the situation is similar to the previous cases (regularity conservation), but for the methods of order $6$ and $8$, resonances appear: for specific values of the stepsize, the regularity of the numerical solution deteriorates. 

\begin{figure}[ht]
\begin{center}
\rotatebox{0}{\resizebox{0.8\linewidth}{0.2\linewidth}{%
   \includegraphics{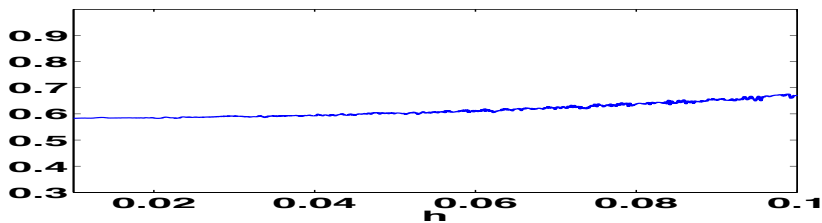}}}
   %\rotatebox{0}{\resizebox{!}{0.3\linewidth}{%
   %\includegraphics{fig8.eps}}}
   \end{center}
   \caption{Order $4$ approximation of the exponential.}
\end{figure}

\begin{figure}[ht]
\begin{center}
\rotatebox{0}{\resizebox{0.8\linewidth}{0.2\linewidth}{%
   \includegraphics{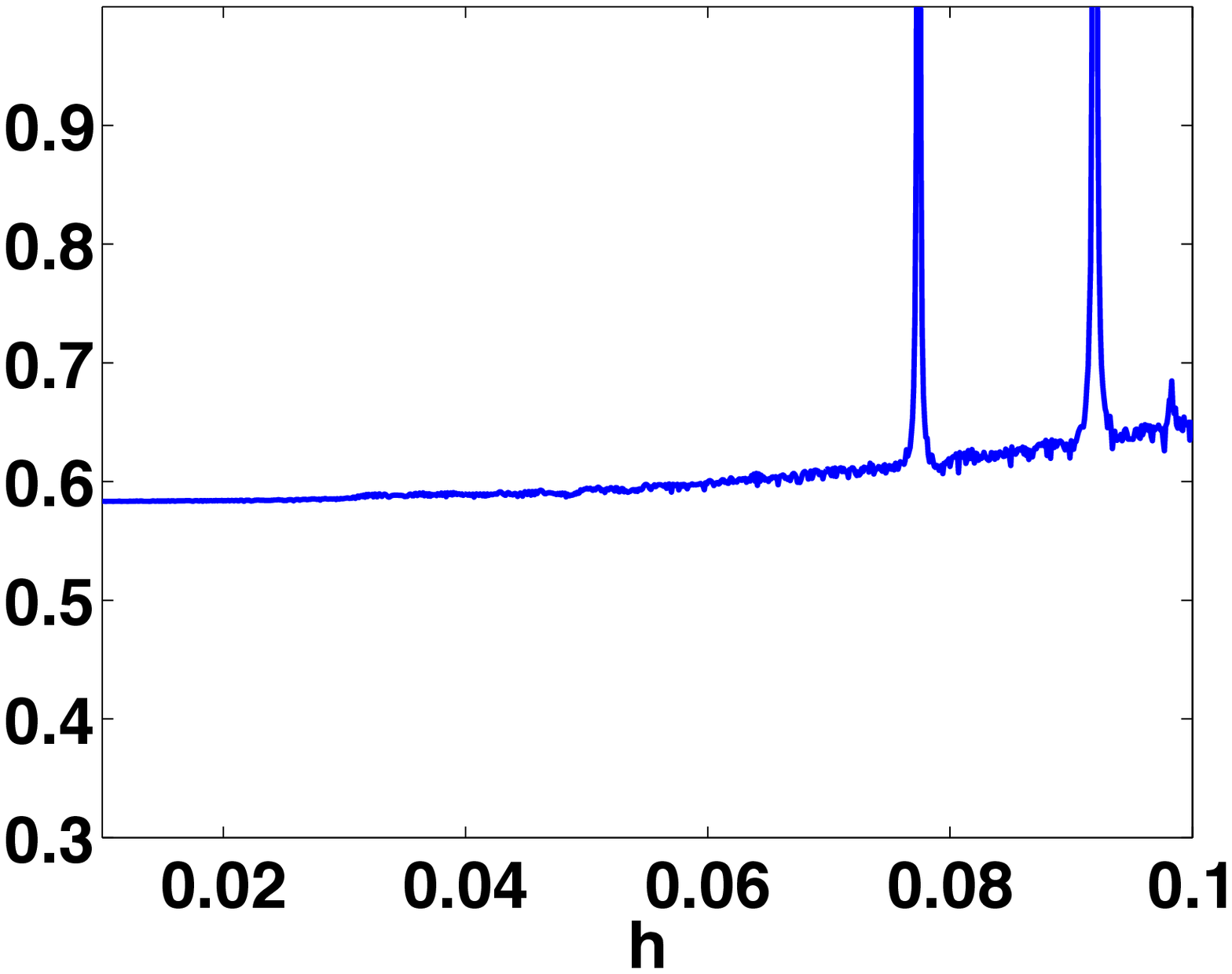}}}
     %\rotatebox{0}{\resizebox{!}{0.3\linewidth}{%
   %\includegraphics{fig8.eps}}}
   \end{center}
   \caption{Order $6$ approximation of the exponential.}
\end{figure}

\begin{figure}[ht]
\begin{center}
\rotatebox{0}{\resizebox{0.8\linewidth}{0.2\linewidth}{%
   \includegraphics{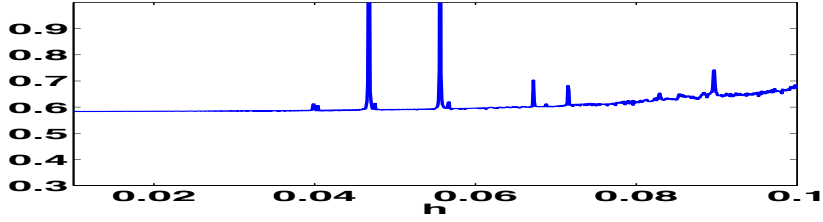}}}
   %\rotatebox{0}{\resizebox{!}{0.3\linewidth}{%
   %\includegraphics{fig8.eps}}}
   \end{center}
   \caption{Order $8$ approximation of the exponential.}
\end{figure}

Finally, we plot in Figure 6 the same simulation for the ``exact'' splitting scheme \eqref{Esplit}. In this last situation, it is known that the resonances appear for step sizes $h$ such that $h (k^2 - \ell^2)$ is close to a multiple of $2\pi$ for some $k$ and $\ell \in \Z$ (see \cite{DF07}). 

\begin{figure}[ht]
\begin{center}
\rotatebox{0}{\resizebox{0.8\linewidth}{0.2\linewidth}{%
   \includegraphics{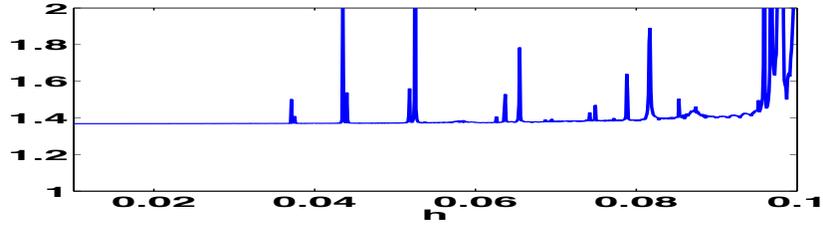}}}
   %\rotatebox{0}{\resizebox{!}{0.3\linewidth}{%
   %\includegraphics{fig8.eps}}}
   \end{center}
   \caption{Exact splitting.}
\end{figure}

The fact that the method of order $4$ possesses a modified energy can easily seen: With the values of $\gamma_1$, $\gamma_2$ and $\gamma_3$ given in \eqref{Egammas}, we have 
$$
R(-\gamma_1 h \Delta)R(-\gamma_2 h \Delta)R(-\gamma_3 h \Delta) = \exp( i  Z_0) 
$$
where
\begin{equation}
\label{Etj}
Z_0 = - 4\arctan \Big(\frac{h\Delta}{2(2 - 2^{1/3})}\Big) + 2 \arctan\Big(\frac{2^{1/3}h\Delta}{2(2 - 2^{1/3})}\Big) = - G(h\Delta/2)
\end{equation}
with
$$
G(x) = 4  \arctan \big(\frac{x}{2 - 2^{1/3}}\big) - 2 \arctan\big(\frac{2^{1/3}x}{2 - 2^{1/3}}\big). 
$$
It is easy to see that for all $x>0$ $G(x)$ is an increasing function such that $G(x) \in (0,\pi)$. Hence Lemma \ref{L1} remains valid for this $Z_0$. Using the same techniques as before, and bounds like 
\eqref{Ebornarc} still valid for the function $G(x)$, we can show the existence of a modified energy for this method, explaining the absence of resonances. 

Note that in the same spirit, we could consider symmetric composition methods based on the order two Strang splitting \eqref{Estrang} to build higher order methods of the form 
$$
\prod_{j = 1}^s
\exp(i\gamma_j hV/2) R(-i\gamma_jh\Delta) \exp(i\gamma_jhV/2)
$$
to approximate \eqref{E1}. A general strategy to show the existence of a modified energy for this method would be to search for an operator $Z(t)$ such that for all $t > 0$, 
$$
\exp(iZ(t)) = \prod_{j = 1}^s
\exp(i\gamma_j tV/2) R(-i\gamma_jh\Delta) \exp(i\gamma_jtV/2)
$$
with 
$$
Z_0 = - \sum_{j =1}^s 2 \arctan(h\gamma_j \Delta/2). 
$$
In the case of the triple jump method, this operator can be written \eqref{Etj}, and the same argument as above shows the existence of a modified energy for this method by using the same kind of techniques. We do not give the details here. The derivation of higher order methods possessing a modified energy is an interesting question that will be addressed in future studies. 

\section*{Acknowledgment}

The authors would like to thank Philippe Chartier for fruitful discussions.  

%\begin{figure}[ht]
%\begin{center}
%\rotatebox{0}{\resizebox{!}{0.3\linewidth}{%
%   \includegraphics{fig4.eps}}}
%   \rotatebox{0}{\resizebox{!}{0.3\linewidth}{%
%   \includegraphics{fig8.eps}}}
%   \end{center}
%   \caption{Maximal oscillations in the truncated $H^1$ norm for $N = 4$ and $N = 8$.}
%\end{figure}

\end{document}